\documentclass[journal,twoside,web]{ieeecolor}
\usepackage{generic}
\usepackage{hyperref}
\usepackage{tabularx}
\usepackage{booktabs}       

\usepackage{textcomp}
\def\BibTeX{{\rm B\kern-.05em{\sc i\kern-.025em b}\kern-.08em
    T\kern-.1667em\lower.7ex\hbox{E}\kern-.125emX}}
\usepackage{amsthm,amsmath,amssymb}

\usepackage{ifpdf}
\usepackage{amsfonts}
\usepackage{mathrsfs}
\usepackage{graphicx}
\usepackage{subfigure}
\usepackage{epstopdf}
\usepackage{lineno}
\usepackage{cite}
\usepackage{enumerate}
\usepackage{algorithm, algorithmic}
\usepackage{bm}
\usepackage{tikz}

\newtheorem{theorem}{Theorem}
\newtheorem{definition}{Definition}

\newtheorem{lemma}{Lemma}     
\newtheorem{remark}{Remark}

\newtheorem{assumption}{Assumption}

\begin{document}
\title{Risk-Averse Learning with Varying Risk Levels}
\author{Siyi Wang, Zifan Wang, Karl H. Johansson, \IEEEmembership{Fellow, IEEE}
\thanks{This work was supported by the Swedish Research Council Distinguished Professor Grant 2017-01078, Knut and Alice Wallenberg Foundation, Wallenberg Scholar Grant, and Swedish Strategic Research Foundation SUCCESS Grant FUS21-0026.}
\thanks{Siyi~Wang, Zifan~Wang and Karl H. Johansson are with the Division of Decision and Control Systems, School of Electrical Engineering and Computer Science, KTH Royal Institute of Technology,
10044  Stockholm, Sweden, e-mail: \{siyiw,zifanw, kallej\}@kth.se.}}
\maketitle

\begin{abstract}
In safety-critical decision-making, the environment may evolve over time, and the learner adjusts its risk level accordingly.
This work investigates risk-averse online optimization in dynamic environments with varying risk levels,  employing Conditional Value-at-Risk (CVaR) as the risk measure.
To capture the dynamics of the environment and risk levels, we employ the function variation metric and introduce a novel risk-level variation metric.
Two information settings are considered: a first-order scenario, where the learner observes both function values and their gradients; and a zeroth-order scenario, where only function evaluations are available.
For both cases, we develop risk-averse learning algorithms with a limited sampling budget and analyze their dynamic regret bounds in terms of function variation, risk-level variation, and the total number of samples.
The regret analysis demonstrates the adaptability of the algorithms in non-stationary and risk-sensitive settings. Finally, numerical experiments are presented to demonstrate the efficacy of the methods.
\end{abstract}

\begin{IEEEkeywords}
Dynamic regret, First-order optimization, Online convex optimization, Risk-averse, 
Zeroth-order optimization
\end{IEEEkeywords}

\section{Introduction}\label{sec:introduction}
Online convex optimization provides a framework for sequential decision-making under uncertainty, where a learner iteratively selects decisions and observes losses revealed by the environment \cite{hazan2016introduction}. 
It finds applications in machine learning \cite{shalev2012online}, finance \cite{li2014online}, and resource allocation \cite{chen2017online}. 
The performance of an optimization algorithm is evaluated by regret, which is defined as the difference between the cumulative loss of the algorithm and that of the best fixed decision in hindsight. In non-stationary environments, dynamic regret, which compares against the sequence of best decisions per iteration, offers a more suitable performance criterion \cite{zinkevich2003online}. 
Additionally, in safety-critical and high-stakes domains such as finance, communication networks, and robotics, extreme losses may lead to catastrophic consequences. Consequently, the primary objective of the learner is to avoid risky outcomes. To manage such risks, common measures include variance, Value-at-Risk (VaR) \cite{larsen2002algorithms}, and Conditional Value-at-Risk (CVaR) \cite{rockafellar2000optimization,kishida2022risk,chapman2022optimizing}. Among these measures, CVaR is widely used because it captures expected tail losses while retaining coherence and convexity for convex losses, thereby ensuring tractable optimization \cite{rockafellar2000optimization}.

In online convex optimization, decisions are typically updated via gradient descent approaches. However, since CVaR depends on the distribution of stochastic costs, its value and exact gradient are typically inaccessible.
This has motivated methods to estimate CVaR gradients from bandit feedback. Existing approaches can be classified into first-order methods \cite{hong2009simulating,kalogerias2022fast,wang2024learning} and zeroth-order methods \cite{cardoso2019risk,wang2024risk}. First-order methods require both function values and gradients, which find applications in online regression \cite{hazan2016introduction} and portfolio selection \cite{helmbold1998line}. Zeroth-order methods assume access only to function evaluations, which find applications in online advertising and dynamic portfolio management \cite{rockafellar2000optimization}. While first-order methods generally achieve faster convergence, they require stronger regularity assumptions \cite{hong2009simulating}.

In sequential decision-making, environments and optimization objectives often evolve dynamically. Such nonstationarity is commonly modeled through time-varying cost functions \cite{besbes2015non}, evolving distributions \cite{cao2020online,pun2023distributionally,wang2024risk}, and dynamic constraints \cite{cao2018online}, with metrics such as function variation and comparator path length to quantify changes. Risk-averse decision-makers are particularly sensitive to such fluctuations. However, only few works address this challenge. For instance, \cite{wang2024risk} considers risk-averse learning with nonstationary distributions for convex and strongly convex objectives, while \cite{pun2023distributionally} studies distributionally time-varying optimization under the Polyak–Łojasiewicz condition.

In dynamic environments, agent may adjust its risk level to ensure consistent performance. For instance, a financial investor may become more conservative in a downturn, or a robot may heighten its risk sensitivity in uncertain terrain. In these scenarios, both the dynamic environment and the adaptive risk level jointly shape the optimization objective. Existing literature on CVaR optimization, such as \cite{borkar2014risk,anderson2020varying,nakagawa2021rm}, typically considers fixed risk levels. A notable exception is \cite{zhu2018simulation}, which qualitatively explores adaptive risk levels to accelerate convergence in CVaR optimization.

In this work, we develop both first- and zeroth-order risk-averse online learning algorithms under varying cost functions and time-varying risk levels.
In the first-order setting, we assume access to multiple function evaluations and their gradients.  To estimate the CVaR gradient, we adopt the dual formulation of CVaR \cite{rockafellar2000optimization}, using the empirical cost distribution built from sampled data and integrating cost gradients to derive the estimator.
Though estimation accuracy increases with sample size, collecting an unlimited number of samples is infeasible. 
Hence, we impose a sampling budget on the total number of samples over the horizon. 
In the zeroth-order setting, only function evaluations are available, which are also subject to sampling budget constraints. A zeroth-order optimization approach is then used to construct the CVaR gradient estimate. 
In both settings, since the empirical distributions are constructed from samples, the CVaR gradient estimation error depends on the discrepancy between the true distribution and its empirical approximation.
Since CVaR variation stems from variations in cost and risk level, we characterize CVaR dynamics using a function variation metric and introduce a new metric to quantify variations in risk level. 
We partition the horizon into intervals, bound the regret in each interval, and aggregate these regret bounds.
By appropriately selecting algorithmic parameters, we achieve a variation-adaptive dynamic regret bound for both the first- and zeroth-order methods.
Our contributions can be summarized as follows: (1) we develop first- and zeroth-order risk-averse learning algorithms that account for both function variation and risk-level variation, which yields a less conservative approach than bounding the overall CVaR variation. These algorithms achieve dynamic regret that adapts to the underlying variation and becomes sublinear when the cumulative variation grows sublinearly, see  Theorems~\ref{theorem:first-order} and \ref{theorem:zeroth-order};
(2) we characterize the performance of the first-order risk-averse algorithm as a function of the number of available samples, revealing how sample size influences the regret.

The remainder of this article is structured as follows: Section~\ref{sec:problem} outlines the preliminaries and problem formulation. Sections~\ref{sec:first-order} and \ref{sec:zeroth-order} present the main results on first- and zeroth-order risk-averse learning with dynamic cost functions and time-varying risk levels. Section~\ref{sec:simulation} demonstrates the efficacy of the algorithms through numerical simulations. Section~\ref{sec:conclusion} concludes this work. Some proofs are given in the Appendix. \\
\textbf{Notations:} Let $\|\cdot\|$ denote the $l_2$ norm.  Let the notation $\mathcal{O}$ hide the constant and $\tilde{\mathcal{O}}$ hide constant and polylogarithmic factors of the number of iterations $T$, respectively.
Let $A \oplus B = \{a + b \vert a \in A, b \in B\}$ denote the Minkowski sum of two sets of position vectors $A$ and $B$ in Euclidean space. 
Let $\mathbb{P}(C)$ denote the probability of the event $C$.

\section{Problem Formulation}
\label{sec:problem}

Consider the cost function $J_t(x,\xi) : \mathcal{X}  \times \Xi 
\rightarrow \mathbb{R} $ for $t = 1, \dots, T$, where  $x\in \mathcal{X}$ denotes the decision variable with $\mathcal{X} \subseteq \mathbb{R}^d$ being the admissible set, $\xi \subseteq  \Xi $ denotes the random noise with the support $\Xi$, and $T$ denotes the horizon. 
Denote the diameter of the admissible set $\mathcal{X} $ as $D_x = \sup_{x,y \in \mathcal{X} } \| x-y\|$. 

To quantify environmental non-stationarity, we employ the function variation metric \cite{besbes2015non} to bound the temporal drift of the cost functions as below.
\begin{definition}\label{def:function variation}
Let $\mathcal{J}$ be the set of cost functions $J_t$, for $t=1,\dots,T$. 
The function variation $V_f$ is defined as drifts of the function sequences $\{J_t\}_{t=1}^T \in \mathcal{J}$  over the  horizon $T$, i.e.,  $  V_f= \sum_{t=2}^T \sup_{x\in\mathcal{X}}\mathbb{E}_\xi\big[|J_t(x,\xi) - J_{t-1}(x,\xi)| \big] $. 
\end{definition} 
\subsection{CVaR}
We use CVaR as risk measure. Suppose $J_t(x,\xi)$ has the cumulative distribution function $F_t(y) = \mathbb{P}(J_t(x,\xi) \le y)$, and is bounded by $U>0$, i.e., $|J_t(x,\xi)|\le U$. Given a confidence level $\alpha \in (0,1]$, the $\alpha$-VaR is
\begin{equation*}
    J_t^{\alpha} = F_t^{-1}(\alpha):= \inf \{y: F_t(y) \geq \alpha \}.
\end{equation*}
The $\alpha$-CVaR describes the expectation of the $\alpha$-fraction of the worst outcomes of $J_t(x,\xi)$ \cite{rockafellar2000optimization}, which is defined as
\begin{align}\label{eq:cvar definition}
\mathrm{CVaR}_{\alpha}\left[J_t(x,\xi)\right]  =\mathbb{E}_{F_t}\left[J_t(x,\xi) \vert J_t(x,\xi) \geq J_t^{\alpha}\right]. 
\end{align}

The following assumptions are widely used in the online convex optimization literature \cite{hazan2016introduction}. 
\begin{assumption}\label{assumption:convex}
The cost function $J_t(x,\xi) $ is convex in $x$ for every $\xi \in \Xi$. 
\end{assumption}
\begin{assumption}\label{assumption:Lipschitz}
The cost function $J_t(x,\xi) $ is Lipschitz continuous in $x$ for every $\xi \in \Xi$. That is, there exists a positive real constant $L_0$ such that, for all $x, y \in \mathcal{X}$, we have $|J_t(x,\xi)-J_t(y,\xi)| \le L_0\|x-y\|$. 
\end{assumption}
The following is a classic result about convexity of the CVaR function \eqref{eq:cvar definition}. 
\begin{lemma}\label{lemma:cvar:convex}\cite{rockafellar2000optimization}
Given Assumption \ref{assumption:convex},
$\mathrm{CVaR}_{\alpha}\left[J_t(x,\xi)\right]$ is convex in $x$.
\end{lemma}
Environmental changes may influence an
agent’s risk tolerance. We consider a varying risk level $\alpha_t$ for $t=1,\dots,T$. Similar to Definition~\ref{def:function variation}, we introduce a risk-level variation metric as below.
\begin{definition}\label{def:risk-level variation} 
Let the risk level be time-varying, denoted as $\alpha_t\in (0,1]$, for $t=1,\dots,T$. The risk-level variation $V_\alpha$ is defined as the drifts of the risk level sequences $\{\alpha_t\}_{t=1}^T $ over the horizon $T$  , i.e.,  $  V_\alpha= \sum_{t=2}^T |\alpha_t - \alpha_{t-1}|.$  
\end{definition}
Another reason for introducing the function and risk-level variation metric is that CVaR is typically not directly observable, whereas cost values (and, in some settings, gradients) and risk levels are. 

Denote $C_t(x)=\mathrm{CVaR}_{\alpha_t}\left[J_t(x,\xi)\right]$. We use the dynamic regret to measure the performance of the algorithm, which is defined as the cumulative loss under the performed actions against the best action per iteration: 
\begin{equation}\label{eq:dynamic regret definition}
    \text{DR}(T)  = \sum_{t=1}^T C_t(x_t)-  \sum_{t=1}^T C_t(x_t^\ast),
\end{equation}
where the action $x_t$ is generated by the algorithm at time $t$, and $ x_t^\ast = \arg\min_{x_t \in \mathcal{X}}C_t (x_t)$ denotes the optimal action at time $t$, for $t=1,\dots, T$. 
Dynamic regret evaluates the algorithm’s ability to track moving optima, 
which is widely used in dynamic settings.

\section{First-order case}\label{sec:first-order}
In this section, we assume access to the stochastic cost $J_t(x,\xi)$ and its gradient, and use these to estimate the CVaR gradient. Then, we develop the first-order risk-averse learning algorithm and analyze its dynamic regret. 

\begin{algorithm}[t] 
\caption{First-order risk-averse learning} \label{alg:first-order}
\begin{algorithmic}[1]
    \REQUIRE Initial value $x_1$, time horizon $T$,  learning rate $\eta$, risk level sequence $\{\alpha_t\}_{t=1}^T$. 
    \FOR{$ {\rm{time}} \;t = 1,\dots, T$} 
    \FOR{$i=1,\ldots,n_t$}
            \STATE Sample $\xi_t^i$, compute cost $J_t(x_{t},\xi_t^i)$ and its gradient $\nabla_{x} J_t(x,\xi_t^i)|_{x=x_t}$
            \ENDFOR
            \STATE Build empirical distribution function $\hat{F}_t^1(y)$, as in \eqref{eq:EDF}
            \STATE Estimate VaR and obtain $ \hat{\nu}_t$, as in \eqref{eq:var estimate} 
            \STATE Construct CVaR gradient estimate $g_t^1$, as in \eqref{eq:first-order gradient estimate} 
            \STATE Update decision: $x_{t+1} \leftarrow \mathcal{P}_{\mathcal{X} } ( x_{t} - \eta  g_t^1)$
    \ENDFOR
\end{algorithmic}
\end{algorithm}

In the following, we introduce the preliminaries for the CVaR gradient estimate. 
The CVaR can be equivalently written as \cite{rockafellar2000optimization}:
\begin{equation}\label{eq:CVaR LP}
{\text{CVaR}}_{\alpha_t}[J_t(x,\xi)] = 
\min_{\nu \in \mathbb{R}}\Big\{ \nu + \frac{1}{\alpha_t}\mathbb{E}_\xi[J_t(x,\xi)- \nu]_{+}\Big\} ,
\end{equation}
where $[x]_+ = \max\{x,0 \}$. 
Denote $\nu_t^\ast(x)$ as the value such that \eqref{eq:CVaR LP} takes the minimum value. It holds that  $\nu_t^\ast(x) = \text{VaR}_{\alpha_t}[J_t(x,\xi)]$, see \cite{rockafellar2000optimization}. 
Define an auxiliary function: 
\begin{equation*} 
    H_t(x,\nu) = \nu + \frac{1}{\alpha_t}\mathbb{E}_\xi[J_t(x,\xi)- \nu]_{+}.
\end{equation*}
From \eqref{eq:CVaR LP}, we have that ${\text{CVaR}}_{\alpha_t}[J_t(x,\xi)] = H_t(x,\nu_t^{\ast}(x))$. Denoting the feasible set $\mathcal{A}_t^\ast(x) := {\text{argmin}}_{\nu}H_t(x,\nu)$, we have $\nu_t^\ast(x) = \text{left~endpoint~of}~\mathcal{A}_t^\ast(x)$ \cite{rockafellar2000optimization}.   

Assume that $J_t(\cdot,\xi)$ is Lipschitz and differentiable for every $\xi$. Then, the gradient of $H_t(x,\nu)$ is given as \cite{kalogerias2022fast}
\begin{equation*} 
    \nabla_{x,\nu} H_t(x,\nu) = \mathbb{E}_{ \xi} \left[ 
    \begin{array}{c}
         \frac{1}{\alpha_t}\mathbf{1}\{J_t(x,\xi) \ge \nu \}\nabla_{x}J_t(x,\xi)  \\
       1-  \frac{1}{\alpha_t}\mathbf{1}\{J_t(x,\xi) \ge \nu \} 
    \end{array}
    \right],
\end{equation*}
where $\nabla_{x,\nu} H_t(x,\nu)$ represents the gradient of $H_t(x,\nu)$ with respect to $(x,\nu)$ and $\nabla_x J_t(x,\xi)$ represents the gradient of $J_t(x,\xi)$ with respect to $x$.  

The following lemma connects the CVaR gradient with the first-order information $\nabla_x J_t(x,\xi)$.
\begin{lemma}\label{lemma:CVaR derivative}
\cite{hong2009simulating}
    It holds that 
\begin{align}\label{eq:CVaR gradient}
\nabla C_t(x) &= \nabla_x H_t(x,\nu)\vert_{\nu=\nu_t^\ast(x)} \nonumber \\  & = \mathbb{E}_{\xi}\bigg[\frac{1}{\alpha_t} \mathbf{1}\{J_t(x,\xi) \ge \nu_t^\ast(x)\} \nabla_x J_t(x,\xi)\bigg],
\end{align}
where $\nabla C_t(x) $ denotes the gradient of $C_t(x)$ with respect to $x$, and $\nabla_x H_t(x,\nu) $ denotes the gradient of $ H_t(x,\nu)$ with respect to $x$.
\end{lemma}

\subsection{First-order algorithm}\label{sec:CVaR gradient estimate}
This section derives CVaR gradient estimates from sampled stochastic costs and their gradients, and then introduces a first-order risk-averse learning algorithm that updates decisions using these estimates.

Assume that the learner accesses $n_t$ random seeds at time $t$, denoted $\xi_t^i$,  $i = 1,\dots,n_t$. 
Assume that the total number of available samples over horizon $T$ satisfies
\begin{align}\label{eq:sampling requirement}
    \sum_{t=1}^{T} \frac{1}{\sqrt{n_t}} \le c  T^{1-\frac{a}{2}}
\end{align}
with the tuning parameter $a>0$ and some constant $c>0$. 
Given the noise samples, we compute the cost function $J_t(x_t, \xi_t^i)$, for $i=1,\dots,n_t$, and  construct the empirical distribution function as
\begin{equation}\label{eq:EDF}
    \hat{F}_t^1(y)=\frac{1}{n_t} \sum_{i=1}^{n_t} \mathbf{1} \{J_t(x_t, \xi_t^i) \leq y \}.
\end{equation} 
Based on \eqref{eq:EDF}, the VaR estimate is given as 
\begin{equation}\label{eq:var estimate}
    \hat{\nu}_t = \text{VaR}_{\alpha_t}[\hat{F}_t^1]. 
\end{equation}
Leveraging Lemma~\ref{lemma:CVaR derivative}, we construct the first-order CVaR gradient estimate as
\begin{equation}\label{eq:first-order gradient estimate}
g^1_t = \frac{1}{n_t\alpha_t }  \sum_{i=1}^{n_t}\textbf{1}\{J_t(x_t,\xi_t^i) \ge \hat{\nu}_t\}\nabla_xJ_t(x,\xi_t^i)|_{x=x_t}.
\end{equation}
To simplify notation, we hereafter denote $\nabla_xJ(x_t,\xi ) := \nabla_xJ(x,\xi)|_{x=x_t}$. 
In \cite{rockafellar2000optimization}, the learning algorithm is assumed to access $t$ samples at time $t$ and constructs the gradient descent as 
\begin{equation*} 
\tilde{g}^1_t = \frac{1}{t\alpha_t } \sum_{k=1}^t \textbf{1}\{J_t(x_t,\xi_k ) \ge \hat{\nu}_t\}\nabla_xJ_t(x_t,\xi_k).
\end{equation*} 
This setting corresponds to a special case of \eqref{eq:first-order gradient estimate} with $n_t=t$ for all 
$t$ and the sampling parameter $a=1$ in \eqref{eq:sampling requirement}.
 
Given the CVaR gradient estimate \eqref{eq:first-order gradient estimate}, gradient descent updates as 
\begin{equation}\label{eq:gradient descent}
    x_{t+1} = \mathcal{P}_{\mathcal{X}}(x_t-\eta g_t^1)
\end{equation}
with initial value $x_1$, where $\mathcal{P}_{\mathcal{X}}(x):= {\text{argmin}}_{y\in \mathcal{X}}\|x - y\|^2$ denotes the projection operator. 
To summarize, the first-order risk-averse learning algorithm updates the decision using the CVaR gradients constructed from sampled costs and gradients. The procedure is outlined in Algorithm~\ref{alg:first-order}.

\subsection{Regret analysis}
This section analyzes the dynamic regret of Algorithm~\ref{alg:first-order}. 

From \eqref{eq:CVaR gradient} and \eqref{eq:first-order gradient estimate}, we have $\mathbb{E}[g_t^1] = \nabla_xH_t(x_t,\hat{\nu}_t)$.  
Thus, the error of the first-order CVaR gradient estimate is given as 
\begin{align}\label{eq:CVaR estimation error 1}
    e_t^1 & =   \nabla_xH_t(x_t,\hat{\nu}_t) -  \nabla_xC_t(x_t) \nonumber \\ 
    & =  \nabla_xH_t(x_t,\hat{\nu}_t) -   \nabla_x H_t(x_t,\nu_t^\ast) \ ,
\end{align}
where $\nu_t^\ast $ is defined after \eqref{eq:CVaR LP}.  
It can be observed from \eqref{eq:CVaR estimation error 1} that the error of the CVaR gradient estimate is determined by the accuracy of the VaR estimate. We provide the following assumptions to bound the error of the VaR estimate.
\begin{assumption}\label{assumption:J gradient bound}
Let $\|\nabla J_t(x,\xi_t)\| \le G$, for all $t=1,\dots,T$, where $G$ is a positive scalar.  
\end{assumption}
\begin{assumption}\label{assumption:probability lower bound}
Let $F_t^x(y) = \mathbb{P}\{J_t(x,\xi) \le y \}$ and $\mathcal{Y}_x = \text{Range}(J_t(x,\xi))$. For every $x \in \mathcal{X}$, the distribution function $F_t^x$ is continuously differentiable and $L_g$-Lipschitz continuous. Moreover,  its probability density function is  lower bounded by $\underline{p}$, i.e., ${F_t^x}'(y) \ge \underline{p}>0$ for all $y \in \mathcal{Y}_x$. 
\end{assumption}
Assumption~\ref{assumption:J gradient bound} bounds the gradient of the stochastic costs, which is widely used in online optimization literature \cite{ba2025doubly}. Assumption~\ref{assumption:probability lower bound} states that the probability distribution function of the stochastic costs is both upper and lower bounded. The same assumption can be found in \cite{zinkevich2003online,hazan2016introduction}.

The following lemma is provided to bound the errors of the VaR estimate and the CVaR gradient estimate. 
\begin{lemma}\label{lemma:var bound}
Let Assumptions \ref{assumption:J gradient bound} and \ref{assumption:probability lower bound} hold. Assume that $n_t$ samples are available at time $t$. Then, it holds that 
\begin{align}\label{eq:DKW 1}
\mathbb{P}\{|\hat{\nu}_t -\nu_t^\ast| >\varepsilon \} \le 2e^{-2n_t\varepsilon^2\underline{p}^2}.
\end{align}
Given a confidence level $\gamma$, for $t=1,\dots,T$, the error of the CVaR gradient estimate \eqref{eq:CVaR estimation error 1} is bounded as 
\begin{align*} 
\|e_t^1\|
 \le \frac{GL_g\sqrt{\ln (2T/\gamma)}}{\alpha_t\underline{p}\sqrt{2n_t} } 
\end{align*}
with probability $1-\gamma$. 
\end{lemma}
  
In Theorem~\ref{theorem:first-order} below, we will show that the regret primarily consists of three components: (i) the CVaR gradient estimate error; (ii) the error due to periodic examination of function and risk-level drifts, and (iii) the gradient descent optimization error. 
To bound the second regret component, we next present two lemmas that quantify the variation of the CVaR functions induced by function drifts and varying risk levels, respectively.
\begin{lemma}\label{lemma:cvar-risk variation bound}
Let  $X$ be a random variable with distribution $\mathcal{D}_x$, and let  $\alpha_1, \alpha_2 \in (0,1]$ be two risk levels. Let $ |X|\leq U$. Then, we have that
\begin{align*}
    \big|{\rm{CVaR}}_{\alpha_1}[X] - {\rm{CVaR}}_{\alpha_2}[X]\big| \leq \left| \frac{1}{\alpha_1} - \frac{1}{\alpha_2} \right| U .
\end{align*}
\end{lemma}
\begin{lemma}\label{lemma:cvar-function variation bound}
Let $J_1(X)$ and $J_2(X)$ be two functions of a random variable $X$ with distribution $\mathcal{D}_x$. Then, we have that 
\begin{align*} 
        &\big|{\mathrm{CVaR}}_\alpha[J_1(X)] -{\mathrm{CVaR}}_\alpha [J_2(X)] \big| \nonumber \\ 
        \le& \frac{1}{\alpha}  \mathop{\mathbb{E}}\limits_{X\sim \mathcal{D}_x}\big[|J_1(X)-J_2(X)|\big].
    \end{align*}
\end{lemma}
Proofs of Lemmas~\ref{lemma:cvar-risk variation bound} and \ref{lemma:cvar-function variation bound} are provided in the Appendix.

Leveraging Lemmas~\ref{lemma:var bound}--\ref{lemma:cvar-function variation bound}, Theorem~\ref{theorem:first-order} analyzes the dynamic regret of Algorithm~\ref{alg:first-order}.
\begin{theorem}\label{theorem:first-order}
Let Assumptions \ref{assumption:convex}--\ref{assumption:probability lower bound} hold.  Suppose that the total number of samples over the horizon $T$ satisfies \eqref{eq:sampling requirement} with a tuning parameter $a>0$. 
By selecting $  \eta = \left(\frac{V_\alpha+V_f}{T}\right)^{\frac{1}{3}}$, Algorithm~\ref{alg:first-order} achieves  ${\rm DR}(T) =  \tilde{\mathcal{O}}(T^{\frac{2}{3}}(V_\alpha+V_f)^{\frac{1}{3}}+ T^{1-\frac{a}{2}}) $ when $a\in(0,\frac{2}{3})$ and ${\rm DR}(T) =  \tilde{\mathcal{O}}(T^{\frac{2}{3}}(V_\alpha+V_f)^{\frac{1}{3}})$ when $a\ge \frac{2}{3}$, with high probability. 
\end{theorem}    
\textit{Proof.} 
To account for the time-varying function and risk level, we partition the horizon $T$ into $s = \left\lceil\frac{T}{\Delta_T}\right\rceil$ intervals of length $\Delta_T \in (1,T)$ and examine the performance of Algorithm~\ref{alg:first-order} on each interval. Each interval is defined as
\begin{align*}
    \mathcal{T}_j  =  \big\{ t : (j-1)\Delta_T < t \le \min\{T, j\Delta_T\} \big\}, \
        {\rm for}~j = 1, \dots, s&,
\end{align*} 
where $j=\left\lceil\frac{t}{\Delta_T}\right\rceil$ denotes the interval  that time $t$ belongs to.
Denote $\tilde{x}_j^\ast$ as the optimal action over the interval $j$, i.e., $\tilde{x}_j^\ast: = \min_{x\in \mathcal{X}}\sum_{t \in \mathcal{T}_j}C_t(x)$. Then, the dynamic regret \eqref{eq:dynamic regret definition} can be written as 
\begin{align}\label{eq:DR 01}
    \mathrm{DR}(T)  &= \sum_{j=1}^s\sum_{t\in \mathcal{T}_j} \Big(\big(C_t(x_t)- C_t(\tilde{x}_j^\ast) \big) + \big( C_t(\tilde{x}_j^\ast) - C_t(x_t^\ast) \big)\Big) \nonumber \\ 
    & = \mathcal{R}_1 +  \mathcal{R}_2
\end{align}
with $\mathcal{R}_1 = \sum_{j=1}^s\sum_{t\in \mathcal{T}_j}  \big(C_t(x_t)- C_t(\tilde{x}_j^\ast) \big)$ and $\mathcal{R}_2 = \sum_{j=1}^s\sum_{t\in \mathcal{T}_j} \big( C_t(\tilde{x}_j^\ast) - C_t(x_t^\ast) \big)$. 
For $\mathcal{R}_1$, we have 
\begin{align}\label{eq:R1}
\mathcal{R}_1 &\le \sum_{j=1}^s\sum_{t\in \mathcal{T}_j} \langle \nabla_xC_t(x_t), x_t-\tilde{x}_j^\ast \rangle \nonumber \\ 
& = \sum_{j=1}^s\sum_{t\in \mathcal{T}_j} \mathbb{E}[\langle g_t^1 -e_t^1,  x_t-\tilde{x}_j^\ast \rangle ],
\end{align}
where the inequality follows from the convexity of $C_t(x)$, as shown in Lemma~\ref{lemma:cvar:convex}, and the equality from substituting \eqref{eq:CVaR estimation error 1} into \eqref{eq:R1}. 
By the gradient descent rule \eqref{eq:gradient descent}, for $t\in \mathcal{T}_j$, we have that 
\begin{flalign*}
 \|x_{t+1} - \tilde{x}_j^{\ast} \|^2  
  =   &\|\mathcal{P}_{\mathcal{X}}(x_{t}-\eta g_t^1) - \tilde{x}_j^{\ast} \|^2  \\
 \leq&  \| x_{t}-\eta g_t^1  - \tilde{x}_j^{\ast}\|^2   \\ 
=&  \|x_t-\tilde{x}_j^{\ast}\|^2+\eta^2\|g_t^1\|^2-2\eta \langle g_t^1, x_t-\tilde{x}_j^{\ast} \rangle,    &&
\end{flalign*}
where the inequality holds as $\tilde{x}_j^{\ast} \in \mathcal{X} $. 
Then, we obtain
\begin{align}\label{eq:gt-x}
  \hspace{-0.5em} \langle g_t^1, x_t-\tilde{x}_j^{\ast} \rangle \le \frac{\|x_t-\tilde{x}_j^{\ast}\|^2 -\|x_{t+1}-\tilde{x}_j^{\ast}\|^2}{2\eta} + \frac{\eta}{2}\|g_t^1\|^2,\hspace{-0.5em}
\end{align}
for $t \in \mathcal{T}_j$.
Substituting \eqref{eq:gt-x}  into $\mathcal{R}_1$, we
obtain 
\begin{align}\label{eq:R1 1}
\mathcal{R}_1   
 &\leq \sum_{j=1}^s\sum_{t\in \mathcal{T}_j} \Big( \frac{1}{2\eta} \mathbb{E}[ \|x_t-\tilde{x}_j^{\ast}\|^2 -\|x_{t+1}-\tilde{x}_j^{\ast}\|^2  ]   \nonumber \\
&\hspace{1em}+ \frac{\eta}{2}\mathbb{E}[\|g_t^1\|^2]  + \mathbb{E}[\| e_t^1\|\| x_t -\tilde{x}_j^{\ast}\| ] \Big)\nonumber  \\
&= \frac{1}{2 \eta } \sum_{j=1}^{s-1}\left( 
 \| x_{(j-1)\Delta_T+1}  -  \tilde{x}_j^{\ast}\|^2 -  \|x_{j\Delta_T} - \tilde{x}_j^{\ast}\|^2 \right) \nonumber \\
 &\hspace{1em}+ \| x_{(s-1)\Delta_T+1}  -  \tilde{x}_s^{\ast}\|^2 -  \|x_{T} - \tilde{x}_s^{\ast}\|^2+ \mathcal{R}_{11} + \mathcal{R}_{12} \nonumber \\
 &\leq \frac{D_x^2T}{ \eta\Delta_T }  + \mathcal{R}_{11} + \mathcal{R}_{12},
\end{align}
with 
$ \mathcal{R}_{11} =\sum_{t=1}^T \frac{\eta}{2}\mathbb{E}[\|g_t^1\|^2] $ and $ \mathcal{R}_{12} = \sum_{j=1}^s \sum_{t\in \mathcal{T}_j}  \mathbb{E}[\| e_t^1\|\| x_t -\tilde{x}_j^{\ast}\| ] $.
The second inequality of \eqref{eq:R1 1} follows from the diameter definition $D_x = \sup_{x,y \in \mathcal{X}} \| x-y\|$ and $\left\lceil\frac{T}{\Delta_T}\right\rceil \le \frac{2T}{\Delta_T} $. By Assumption~\ref{assumption:J gradient bound}, we have $\mathcal{R}_{11} \le \frac{\eta G^2T}{2}$.
For $\mathcal{R}_{12}$, we have 
\begin{align}\label{eq:R12}
 \mathcal{R}_{12} &\le D_x    \sum_{t=1}^T \mathbb{E}[\|e_t^1\|]   \le \frac{D_xGL_g\sqrt{\ln (2T/\gamma)}}{\sqrt{2}\alpha_{\min} \underline{p} } \sum_{t=1}^T \frac{1}{\sqrt{n_t}}  \nonumber \\ 
 &\le  \frac{c D_x G L_g\sqrt{\ln{(2T/\gamma)}}}{\sqrt{2}\alpha_{\min} \underline{p} }  T^{1-\frac{a}{2}}
\end{align}
with probability $1-\gamma$. In \eqref{eq:R12}, the first inequality follows from $D_x = \sup_{x,y \in \mathcal{X}} \| x-y\|$, the second inequality from Lemma~\ref{lemma:var bound}, and the last inequality from the sampling requirement \eqref{eq:sampling requirement}. 

We next bound $\mathcal{R}_2$. For $t=1,\dots, T$, replacing the function $f_t$ of Lemma~\ref{lemma:function variation} (in the Appendix) with $C_t$, we obtain 
\begin{align}\label{eq:R2}
 \mathcal{R}_{2}    & \le  2  \Delta_T\sum_{t=2}^T \sup_{x\in \mathcal{X}}|C_t(x)-C_{t-1}(x)| \nonumber \\
& =2  \Delta_T\sum_{t=2}^T\sup_{x\in \mathcal{X}} \bigg| \text{CVaR}_{\alpha_t}[J_t(x)] - \text{CVaR}_{\alpha_{t-1}}[J_{t-1}(x)] \bigg|\nonumber \\ 
& =2  \Delta_T\sum_{t=2}^T\sup_{x\in \mathcal{X}} \bigg| \text{CVaR}_{\alpha_t}[J_t(x)] - \text{CVaR}_{\alpha_{t-1}}[J_t(x)] 
\nonumber \\ 
&\hspace{1.5em}
+ \text{CVaR}_{\alpha_{t-1}}[J_t(x)] - \text{CVaR}_{\alpha_{t-1}}[J_{t-1}(x)] \bigg|
\nonumber \\ 
& \le \mathcal{R}_{21} + \mathcal{R}_{22},   
\end{align}
where 
\begin{align*}
\mathcal{R}_{21} &= 2\Delta_T\sum_{t=2}^T\sup_{x\in \mathcal{X}} \big| \text{CVaR}_{\alpha_t}[J_t(x)] - \text{CVaR}_{\alpha_{t-1}}[J_t(x)] \big|,    \nonumber \\ 
\mathcal{R}_{22} &= 2\Delta_T\sum_{t=2}^T\sup_{x\in \mathcal{X}} \big| \text{CVaR}_{\alpha_{t-1}}[J_t(x)] - \text{CVaR}_{\alpha_{t-1}}[J_{t-1}(x)] \big|.
\end{align*}
By Lemma~\ref{lemma:cvar-risk variation bound}, we have 
\begin{equation}\label{eq:R21}
\mathcal{R}_{21} \le 2 \Delta_T  \sum_{t=2}^T  U \bigg|\frac{1}{\alpha_t} -\frac{1}{\alpha_{t-1}}  \bigg|
\le \frac{2 U \Delta_T V_\alpha}{\alpha_{\min}^2} ,
\end{equation}
where $\alpha_{\min} := \min\{a_1,\dots,\alpha_T \}$, and the second inequality holds by Definition~\ref{def:risk-level variation}. 
Furthermore, by Lemma~\ref{lemma:cvar-function variation bound}, we have 
\begin{align}\label{eq:R22}
\mathcal{R}_{22} &\le   \sum_{t=2}^T\frac{2 \Delta_T}{\alpha_{t-1}}\sup_{x}\mathbb{E}_{\xi}\big[|J_t(x,\xi)-J_{t-1}(x,\xi)| \big]  \le  \frac{2 \Delta_TV_f}{\alpha_{\min}},
\end{align}
where the second inequality holds by Definition~\ref{def:function variation}. 
Substituting \eqref{eq:R1 1}--\eqref{eq:R22} into \eqref{eq:DR 01}, we obtain 
\begin{align}\label{eq:DR bound first-order}
 \mathrm{DR}(T) \le& \frac{D_x^2T}{ \eta\Delta_T }+\frac{\eta G^2T}{2}+ \frac{2 U \Delta_T V_\alpha}{\alpha_{\min}^2}  +  \frac{2 \Delta_TV_f}{\alpha_{\min}}\nonumber \\ 
 &+\frac{c D_x G L_g\sqrt{\ln{(2T/\gamma)}}}{\sqrt{2}\alpha_{\min} \underline{p} }  T^{1-\frac{a}{2}}
\end{align}
with probability at least $ 1- \gamma $.
We select $  \eta = \left(\frac{V_\alpha+V_f}{T}\right)^{\frac{1}{3}}$ and $\Delta_T = \left(\frac{T}{V_\alpha+V_f}\right)^{\frac{2}{3}} $ to minimize the regret bound. Then, ${\rm DR}(T) =  \tilde{\mathcal{O}}(T^{\frac{2}{3}}(V_\alpha+V_f)^{\frac{1}{3}}, T^{1-\frac{a}{2}}) $, which reduces to ${\rm DR}(T) =  \tilde{\mathcal{O}}(T^{\frac{2}{3}}(V_\alpha+V_f)^{\frac{1}{3}})$  when $a\ge \frac{2}{3}$, with probability at least $ 1- \gamma $.  \hfill $\qed$
\begin{remark}
When the risk level is static, i.e., $\alpha_t = \alpha$ for all $t = 1,\dots, T$, we have $V_\alpha = 0$. Additionally, when ${\alpha_t}$ is monotonically nonincreasing or nondecreasing with time, we have $V_\alpha = |\alpha_1- \alpha_T|$, which is constant and independent of $T$. In these two cases, Algorithm~\ref{alg:first-order} achieves ${\rm DR}(T) =  \tilde{\mathcal{O}}(T^{\frac{2}{3}}V_f^{\frac{1}{3}}+ T^{1-\frac{a}{2}}) $ when $a\in(0,\frac{2}{3})$ and ${\rm DR}(T) =  \tilde{\mathcal{O}}(T^{\frac{2}{3}}V_f^{\frac{1}{3}})$  when $a\ge \frac{2}{3}$, with probability at least $ 1- \gamma $. 
This argument also applies to a static function setting.
Namely, 
when the function is static, i.e., $J_t(\cdot) = J(\cdot)$ for all $t = 1,\dots, T$, we have $V_f = 0$. Then, Algorithm~\ref{alg:first-order} achieves ${\rm DR}(T) =  \tilde{\mathcal{O}}(T^{\frac{2}{3}}V_\alpha^{\frac{1}{3}}+ T^{1-\frac{a}{2}}) $ when $a\in(0,\frac{2}{3})$ and ${\rm DR}(T) =  \tilde{\mathcal{O}}(T^{\frac{2}{3}}V_\alpha^{\frac{1}{3}})$  when $a\ge \frac{2}{3}$, with probability at least $ 1- \gamma $.    
\end{remark}
\begin{remark}\label{remark:tighter bound}
We characterize the dynamic regret in terms of the variation in the cost function and the variation in the risk level, rather than the variation in CVaR. As noted, the CVaR value is generally not directly observable and has to be estimated from samples. Moreover, variations in both the cost function and risk level will lead to variations in CVaR. In  Theorem~\ref{theorem:first-order}, $V_f$ and $V_\alpha$ influence the regret bound independently, and the dominant-order term governs the asymptotic rate. Consequently, defining the performance metric using these two variations yields a tighter regret bound than that based on the variation in CVaR. 
\end{remark}
We provide an example to further illustrate Remark~\ref{remark:tighter bound}. Let the function vary for $t \in [1,T/2]$ and remain constant for $t \in (T/2, T]$, with the function variation $V_f = \tilde{O}(T^{\frac{1}{2}})$. Let the risk level remain constant for $t \in [1,T/2]$ and vary for $t \in (T/2, T]$, with the variation $V_\alpha = \tilde{O}(T^{\frac{1}{3}})$. 
By Theorem~\ref{theorem:first-order}, the resulting regret bound is dominated by the function variation $V_f$. When regret is expressed directly in terms of CVaR variation, the bound increases, as both function drift and risk-level changes induce variation in the CVaR values.

 \begin{remark}
A natural question arises: \textit{Can the optimal decisions remain unchanged as the CVaR risk level $\alpha_t$ varies}? Characterizing such regimes is nontrivial because CVaR depends on the loss quantile, i.e., the inverse cumulative distribution function, which rarely admits a closed-form expression. This dependence makes sensitivity analysis with respect to $\alpha_t$ analytically challenging. Nonetheless, some tractable special cases exist. For instance, when the random component is separable from the decision variable, i.e., when $J(x, \xi) = f(x) + g(\xi)$, the optimal decision is unaffected by changes in $\alpha_t$. Similarly, in constrained dynamic portfolio problems \cite{rockafellar2000optimization} with linear transaction costs and Gaussian returns, varying $\alpha_t$ merely rescales the risk term, rendering the optimizer invariant across different risk levels.
\end{remark}
\begin{remark}
Theorem~\ref{theorem:first-order} reveals that the regret bound decreases as the total number of samplings (controlled by the parameter $a$) increases when $a \in (0,\frac{2}{3}]$, and ceases to improve when $a > \frac{2}{3}$. This occurs because the VaR estimate error becomes negligible compared to other error terms in \eqref{eq:DR bound first-order} for $a > \frac{2}{3}$. 
Furthermore, the sampling strategy in \cite{rockafellar2000optimization} corresponds to the special case $a=1$ in \eqref{eq:sampling requirement}, which
 yields a convergence rate of $\tilde{\mathcal{O}}(T^{-1/2})$ for a static problem. Under identical static conditions, Algorithm~\ref{alg:first-order} achieves a regret bound of $\mathcal{R}(T) = \tilde{\mathcal{O}}(T^{1/2} + T^{1 - a/2})$, which simplifies to $\tilde{\mathcal{O}}(T^{1/2})$ when $a \ge 2/3$. Consequently, our algorithm matches the $\tilde{\mathcal{O}}(T^{1/2})$ bound of \cite{rockafellar2000optimization} with fewer samples.
\end{remark}

\section{Zeroth-order case}\label{sec:zeroth-order}
This section addresses a bandit setting, where only the function evaluations of the selected action can be observed. We tackle this limitation by adopting a zeroth-order optimization approach to estimate the CVaR gradient.

\begin{algorithm}[t] 
\caption{Zeroth-order risk-averse learning} \label{alg:zeroth-order}
\begin{algorithmic}[1]
    \REQUIRE Initial value $x_1$, time horizon $T$, smoothing parameter $\delta$, learning rate $\eta$, risk level sequence $\{\alpha_t\}_{t=1}^T$.
    \FOR{$ {\rm{time}} \;t = 1,\dots, T$} 
    \STATE  Sample $u_{t} \in \mathbb{S}^{d}$
    \STATE  Perturb $\hat{x}_{t}=x_{t}+\delta u_{t} $
    \FOR{$i=1,\ldots,n_t$}
    \STATE Play $\hat{x}_{t}$ and obtain $J_t(\hat{x}_{t},\xi_t^i)$
    \ENDFOR
    \STATE Build empirical distribution function $\hat{F}_t^0(y)$, as in \eqref{eq:EDF 0}
    \STATE Estimate CVaR: $ {\rm{CVaR}}_{\alpha}[\hat{F}_t^0] $ 
    \STATE Construct gradient estimate $g_t^0$, as in \eqref{eq:gradient descent 0}
    \STATE Update $x$: $x_{t+1} \leftarrow \mathcal{P}_{\mathcal{X}^{\delta}} ( x_{t} - \eta g_t^0)$
    \ENDFOR
\end{algorithmic}
\end{algorithm}

To begin with, we construct a smoothed approximation of the CVaR function \cite{cardoso2019risk}. Assume that the feasible set $\mathcal{X} $ contains the ball of radius $r$ centered at the origin, i.e., $r \mathbb{B} \subseteq  \mathcal{X} $.  
For a given point $x\in \mathcal{X}$, define the perturbed action 
$\hat{x}=x+\delta u$, where $u$ is a direction vector sampled from a unit sphere $\mathbb{S}^d \in \mathbb{R}^d$ and $\delta>0$ is the perturbation radius, which is also known as the smoothing parameter. Then, the smoothed approximation of the CVaR function is given as 
\begin{equation}\label{eq:smoothed cvar}
    C_t^\delta(x) = \mathbb{E}_{u \sim \mathbb{S}^d}[C_t(x+\delta u)].
\end{equation} 
The following lemma characterizes properties of the smoothed CVaR approximation.
\begin{lemma}\label{lemma:smoothed cvar} \cite{cardoso2019risk} Given Assumptions~\ref{assumption:convex} and \ref{assumption:Lipschitz}, we have that  
\begin{enumerate}
    \item $C_t^\delta\left(x\right)$ is convex in $x$;
    \item
    $C_t^\delta (x)$ is $L_0$-Lipschitz in $x$ and $|C_t^\delta(x) - C_t(x)| \le \delta L_0 $.
\end{enumerate}
\end{lemma}

In the zeroth-order setting, the algorithm actually executes the perturbed action $\hat{x}_t$ at time $t$, rather than the nominal center $x_t$. As a consequence, it only observes the cost $C_t(\hat{x}_t)$ instead of $C_t(x_t)$.  Accordingly, we modify the definition of dynamic regret as
\begin{equation}\label{eq:dynamic regret definition 0}
    \text{DR}_0(T)  = \sum_{t=1}^T C_t(\hat{x}_t)-  \sum_{t=1}^T C_t(x_t^\ast),
\end{equation}
where the action $\hat{x}_t$ is generated by the algorithm at time $t$.   
\subsection{Zeroth-order algorithm}
At each time $t$, we execute the perturbed action $\hat{x}_t = x_t +\delta u_t$ for $n_t$ times and obtain the cost  samples $J(\hat{x}_t,\xi_t^i)$, for $i=1,\ldots,n_t$.
Then, we use the queried function values to construct the empirical distribution function
\begin{equation}\label{eq:EDF 0}
    \hat{F}_t^0(y)=\frac{1}{n_t} \sum_{i=1}^{n_t} \mathbf{1} \{J_t(\hat{x}_t, \xi_t^i) \leq y \}.
\end{equation} 
Given \eqref{eq:EDF 0}, we compute the CVaR estimate $\text{CVaR}_{\alpha_t}[\hat{F}_t^0]$. Leveraging the zeroth-order gradient identity  $\nabla f(x) = \mathbb{E}[f(x+\delta u)u]d/\delta$ \cite{flaxman2004online}, the CVaR gradient estimate is constructed as
\begin{equation}\label{eq:estimated gradient 0}
    g_t^0=\frac{d}{\delta} \text{CVaR}_{\alpha_t}\big[\hat{F}_t^0\big] u_t.
\end{equation}
The gradient descent for the risk-averse learning proceeds as 
\begin{equation}\label{eq:gradient descent 0}
    x_{t+1} = \mathcal{P}_{\mathcal{X}^\delta}(x_t - \eta g_t^0),
\end{equation} 
with the initial value $x_1$, where $\mathcal{P}_{\mathcal{X}^\delta }(x):= {\text{argmin}}_{y\in \mathcal{X}^\delta}\|x - y\|^2$ denotes the projection operator with $\mathcal{X}^\delta = \{x \in \mathcal{X} \vert \frac{1}{1-\delta/r} x \in \mathcal{X}\}$ being the  projection set. The projection ensures that the sampled actions $\hat{x}_t$ remain within the admissible set $\mathcal{X}$, which establishes as $\big(1-\frac{\delta}{r}\big)\mathcal{X} \oplus \delta   \mathbb{B}  = \big(1-\frac{\delta}{r}\big)\mathcal{X} \oplus \frac{\delta}{r} r   \mathbb{B}  \subseteq \big(1-\frac{\delta}{r}\big)\mathcal{X} \oplus \frac{\delta}{r} \mathcal{X} = \mathcal{X}.$
The zeroth-order risk-averse learning algorithm updates decisions using the CVaR gradient estimate constructed from stochastic costs, which is summarized in  Algorithm~\ref{alg:zeroth-order}

\subsection{Regret analysis}
The empirical distribution function \eqref{eq:EDF 0} is constructed using stochastic costs, which induces the CVaR gradient estimate error: 
\begin{equation}\label{eq:CVaR estimation error 0}
    e_t^0:= \text{CVaR}_{\alpha_t}[\hat{F}_t^0] - \text{CVaR}_{\alpha_t}[F_t].
\end{equation}
The following lemma bounds \eqref{eq:CVaR estimation error 0}.
\begin{lemma}\label{lemma:cvar-estimation error bound} \cite{wang2022zeroth} Let $F$ and $G$ be the cumulative distribution functions of two random variables, respectively. Let the random variables be bounded by $U$. Then, we have that 
\begin{equation*} 
        |{\mathrm{CVaR}}_\alpha[F] -{\mathrm{CVaR}}_\alpha [G] | \le \frac{U}{\alpha} \sup_{x} |F(x)-G(x)|.
    \end{equation*}
\end{lemma}
Using the preceding results, Theorem~\ref{theorem:zeroth-order} analyzes the dynamic regret of Algorithm~\ref{alg:zeroth-order}.
\begin{theorem}\label{theorem:zeroth-order}
Let Assumptions \ref{assumption:convex} and \ref{assumption:Lipschitz} hold.  Suppose that the total sampling number over the horizon $T$ satisfy \eqref{eq:sampling requirement} with a constant $a>0$. 
\begin{enumerate}
    \item When $a \in (0,\frac{4}{5}]$, by selecting $  \delta =  T^{-\frac{a}{4}}(V_\alpha+V_f)^{\frac{1}{5}}$ and $\eta =T^{-\frac{3a}{4}}(V_\alpha+V_f)^{\frac{3}{5}}$,   
Algorithm~\ref{alg:zeroth-order}  achieves $\rm{DR}_0(T) =  \tilde{\mathcal{O}}(T^{1-\frac{a}{4}}(V_\alpha+V_f)^{\frac{1}{5}})$ with high probability.
 \item When $a >\frac{4}{5} $, by selecting  $  \delta =  T^{-\frac{1}{5}}(V_\alpha+V_f)^{\frac{1}{5}}$ and $\eta =T^{-\frac{3}{5}}(V_\alpha+V_f)^{\frac{3}{5}}$, Algorithm~\ref{alg:zeroth-order}  achieves $\rm{DR}_0(T) =  \tilde{\mathcal{O}}(T^{\frac{4}{5}}(V_\alpha+V_f)^{\frac{1}{5}}) $ with high probability.
\end{enumerate}
\end{theorem}
\textit{Proof.}
The dynamic regret \eqref{eq:dynamic regret definition 0} can be further written as 
\begin{align}\label{eq:convex 11}
\hspace{-0.5em}        \mathrm{DR}_0(T)  & =  \sum_{t=1}^T \mathbb{E}_{u}[C_t(\hat{x}_t+\delta u-\delta u)] - \mathbb{E}_{u}[C_t(x_t^\ast+\delta u-\delta u)]  \nonumber \\
&\leq \sum_{t=1}^T \Big(C_t^\delta( \hat{x}_t)-   C_t^\delta (x_t^\ast) \Big)+2 \delta L_0 T \nonumber \\ 
        &\leq \sum_{t=1}^T \Big(C_t^\delta( x_t)-  C_t^\delta (x_t^\ast)\Big)+3 \delta L_0 T,
\end{align}
where the first inequality follows from the Lipschitzness  of $C_t$ and the definition of $C_t^\delta$, as in \eqref{eq:smoothed cvar}. The second inequality follows from the Lipschitzness of $C_t^\delta$, as shown in Lemma \ref{lemma:smoothed cvar}. 
Furthermore, for $t=1,\dots,T$, we have 
\begin{align}\label{eq:projection mismatch}
   \min_{x_t \in \mathcal{X}^\delta} C_t^\delta(x_t) &=    \min_{x_t \in \mathcal{X}} C_t^\delta\big((1-\delta/r)x_t\big) \nonumber \\
   &\leq  \min_{x_t \in \mathcal{X}} (\delta/r) C_t^\delta(0) + (1-\delta/r)C_t^\delta(x_t) \nonumber \\
    &\leq \min_{x_t\in \mathcal{X}} C_t^\delta(x_t)+ (\delta/r)L_0 \left\|x_t \right\| \nonumber \\
   &\leq \min_{x_t\in \mathcal{X}} C_t^\delta(x_t)+D_x L_0 \delta/r,
\end{align}
where the first inequality follows from the convexity of $C_t^\delta$ and the second inequality from Lipschitzness of $C_t^\delta$, as shown in Lemma~\ref{lemma:smoothed cvar}.  Denote $x_t^{\delta,\ast} = {\textrm{argmin}}_{x \in \mathcal{X}^\delta} C_t^\delta (x)$ as the single-step optimal decision, for $t=1,\dots,T$. 
Substituting \eqref{eq:projection mismatch} into \eqref{eq:convex 11} and following the derivation of \eqref{eq:DR 01}, we have that 
\begin{align}\label{eq:convex 12}
\mathrm{DR}_0(T)  &\leq \sum_{t=1}^T C_t^\delta(x_t)-  \sum_{t=1}^T C_t^\delta (x_t^{\delta,\ast})+(3 +D_x/r)\delta L_0 T    \nonumber \\
&= \tilde{\mathcal{R}}_1 + \tilde{\mathcal{R}}_2 +(3 +D_x/r)\delta L_0 T  
\end{align}
with 
\begin{align*}
 \tilde{\mathcal{R}}_1 &= \sum_{j=1}^s\sum_{t\in \mathcal{T}_j}  \big(C_t^\delta(x_t)- C_t^\delta(\bar{\tilde{x}}_j^{\delta,\ast}) \big) \\
\tilde{\mathcal{R}}_2 &= \sum_{j=1}^s\sum_{t\in \mathcal{T}_j} \big( C_t^\delta(\bar{\tilde{x}}_j^{\delta,\ast}) - C_t^\delta(x_t^{\delta,\ast}) \big) 
\end{align*} and $\bar{\tilde{x}}_j^{\delta,\ast}=\min_{x\in\mathcal{X}^\delta}\sum_{t\in \mathcal{T}_j}C_t^\delta(x)$ is the optimal action over the interval $j$. 
From \eqref{eq:estimated gradient 0} and \eqref{eq:CVaR estimation error 0}, we have that \begin{equation}\label{eq:cvar derivative}
   \nabla C_t^\delta(x_t)   = \mathbb{E}\big[g_t^0 - \frac{d}{\delta}e_t^0u_t\big].
\end{equation}
Furthermore, for $\tilde{R}_1$,  we have that   
\begin{align}\label{eq:R1 zero}
\tilde{R}_1 &\le \sum_{j=1}^s\sum_{t\in \mathcal{T}_j}\left\langle \nabla C^\delta_t(x_t),x_t-\bar{\tilde{x}}_t^{\delta,\ast}   \right\rangle  \nonumber \\ 
& = \sum_{j=1}^s\sum_{t\in \mathcal{T}_j} \mathbb{E}\left[ \big\langle g_t^0 -\frac{d}{\delta}e_t^0u_t,x_t-\bar{\tilde{x}}_t^{\delta,\ast}   \big\rangle \right]
\nonumber \\ 
&\le \frac{D_x^2T}{ \eta\Delta_T }  + \tilde{\mathcal{R}}_{11} + \tilde{\mathcal{R}}_{12},    
\end{align}
with  $\tilde{\mathcal{R}}_{11} = \sum_{t=1}^T \frac{\eta}{2}\mathbb{E}[\|g_t^0\|^2]$, $    \tilde{\mathcal{R}}_{12} = \sum_{j=1}^s \sum_{t\in \mathcal{T}_j} \frac{d}{\delta}  \mathbb{E}[\| e_t^0\|\| x_t -\bar{\tilde{x}}_j^{\delta,\ast}\| ],$ where the first inequality follows from the convexity of $C_t^\delta$, the equality from substituting \eqref{eq:cvar derivative} into \eqref{eq:R1 zero}, and the second inequality follows the derivation of \eqref{eq:R1}-\eqref{eq:R1 1}. 
By the definition of the CVaR gradient estimate \eqref{eq:estimated gradient 0}, we have 
\begin{equation}\label{eq:R11 0}
 \tilde{\mathcal{R}}_{11} \le  \sum_{t=1}^T \frac{\eta}{2}  \left(\frac{d U}{\delta}\right)^2 = \frac{\eta d^2 U^2T}{2\delta^2}  .  
\end{equation}
We next bound $\tilde{\mathcal{R}}_{12}$. 
By the Dvoretzky–Kiefer–Wolfowitz (DKW) inequality \cite{dvoretzky1956asymptotic}, we have that 
\begin{equation}\label{eq:DKW 0}
    \mathbb{P}\left\{ \sup_{y} |\hat{F}_t^0(y) - F_t(y)| \ge \sqrt{\frac{\ln (2 / \bar{\gamma})}{2 n_t}}  \right\} \le  \bar{\gamma}.
\end{equation}
Denote the event in \eqref{eq:DKW 0} as $B_t$, and $\mathbb{P}\{B_t\}$ the occurrence probability of $B_t$, for $t=1,\dots, T$.
By  Lemma~\ref{lemma:cvar-estimation error bound} and the definition of $e_t^0$, as in \eqref{eq:CVaR estimation error 0}, we have 
\begin{align}
\label{eq:DKW cvar mismatch 0}
     |e_t^0 |  \le \frac{U}{\alpha_t}\sup_y\big| \hat{F}_t^0(y) - F_t(y)\big|  \le \frac{U}{\alpha_t} \sqrt{\frac{\ln (2 / \bar{\gamma})}{2 n_t}} 
\end{align}
with probability at least $1 -\bar{\gamma}$, for $t = 1,\dots, T$.  Let $\gamma = \bar{\gamma}T$. 
Then, subsituting \eqref{eq:DKW cvar mismatch 0} into $\tilde{\mathcal{R}}_{12}$, we obtain 
\begin{align}\label{eq:R12 0}
 \tilde{\mathcal{R}}_{12} &    \le  \sum_{t=1}^{T}\frac{dUD_x}{\delta \alpha_t} \sqrt{\frac{\ln (2T / \gamma)}{2 n_t}}    \le \frac{cdUD_x\sqrt{\ln (2T /  \gamma)}}{\sqrt{2}\delta \alpha_{\min}}  T^{1-\frac{a}{2}}.
\end{align}
with probability $1-\gamma$, which establishes as $    1 - \mathbb{P}\{\bigcup_{t=1}^{T} B_t \} \ge 1 - \sum_{t=1}^{T} \mathbb{P}\{B_t\} \ge 1 - T\frac{\gamma}{T} \ge 1-\gamma$.  The second inequality of \eqref{eq:R12 0} follows from the sampling requirement \eqref{eq:sampling requirement}. 
Furthermore, following the derivation of $\mathcal{R}_2$, as in \eqref{eq:R2}--\eqref{eq:R22}, we have 
\begin{equation}\label{eq:R2 0}
\tilde{R}_2 \le    \frac{2 U \Delta_T V_\alpha}{\alpha_{\min}^2}  +  \frac{2 \Delta_TV_f}{\alpha_{\min}}. 
\end{equation} 
Substituting \eqref{eq:R1 zero},  \eqref{eq:R11 0}, \eqref{eq:R12 0}, and \eqref{eq:R2 0} into \eqref{eq:convex 12}, we obtain 
\begin{align*}
\text{DR}_0(T) &\le (3 +D_x/r)\delta L_0 T + \frac{D_x^2T}{ \eta\Delta_T} + \frac{\eta d^2 U^2T}{2\delta^2} \nonumber \\ 
&\hspace{-3em}+ \frac{cdUD_x\sqrt{\ln (2T / \gamma)}}{\sqrt{2}\delta \alpha_{\min}}  T^{1-\frac{a}{2}} +  \frac{2 U \Delta_T V_\alpha}{\alpha_{\min}^2}+ \frac{2 \Delta_TV_f}{\alpha_{\min}}
\end{align*}
with probability $1-\gamma$. 
To minimize the regret $\text{DR}_0(T)$, we select $  \delta =  T^{-\frac{a}{4}}(V_\alpha+V_f)^{\frac{1}{5}}$, $\eta =T^{-\frac{3a}{4}}(V_\alpha+V_f)^{\frac{3}{5}}$, and $\Delta_T = T^{a}(V_\alpha+V_f)^{-\frac{4}{5}}$ when $a \in (0,\frac{4}{5}]$ and $  \delta =  T^{-\frac{1}{5}}(V_\alpha+V_f)^{\frac{1}{5}}$, $\eta =T^{-\frac{3}{5}}(V_\alpha+V_f)^{\frac{3}{5}}$, and $\Delta_T = T^{\frac{4}{5}}(V_\alpha+V_f)^{-\frac{4}{5}}$ when $a >\frac{4}{5} $. 
Hence,
$\text{DR}_0(T) =  \tilde{\mathcal{O}}(T^{1-\frac{a}{4}}(V_\alpha+V_f)^{\frac{1}{5}}) $ when $a\in(0,\frac{4}{5}]$ and $\text{DR}_0(T) =  \tilde{\mathcal{O}}(T^{\frac{4}{5}}(V_\alpha+V_f)^{\frac{1}{5}}) $ when $a > \frac{4}{5}$, with probability at least $ 1- \gamma $.
\hfill $\qed$

Unlike \cite{wang2024risk}, which periodically resets parameters such as the learning rate to handle nonstationarity, we set the learning rate and smoothing parameter as constants. Instead, we partition the time horizon for the purpose of regret analysis only. In this case, by enforcing the sampling requirement over the entire horizon $T$ rather than on each interval $\mathcal{T}_j$, we obtain tighter estimation error bounds (see \eqref{eq:R12} and \eqref{eq:R12 0}) and, consequently, a sharper overall regret bound. 

\begin{table}[h]
\label{table:1}
\centering{\caption{Order of Regret}}
\begin{tabularx}{0.48\textwidth}{ 
   >{\centering\arraybackslash}X 
    >{\centering\arraybackslash}X 
   >{\centering\arraybackslash}X  } 
\toprule[1pt]
 Sampling parameter & First-order  &Zeroth-order\\
\midrule[1pt]  
$0<a\le\frac{2}{3} $ &  $T^{\frac{2}{3}}V_T^{\frac{1}{3}}+ T^{1-\frac{a}{2}}$   & $T^{1-\frac{a}{4}}V_T^{\frac{1}{5}}$ \\ 
\hline 
$\frac{2}{3} <a \le \frac{4}{5}$ &  $T^{\frac{2}{3}}V_T^{\frac{1}{3}}$   &  $T^{1-\frac{a}{4}}V_T^{\frac{1}{5}}$  \\ 
\hline 
$a>\frac{4}{5} $ & $T^{\frac{2}{3}}V_T^{\frac{1}{3}}$   & $T^{\frac{4}{5}}V_T^{\frac{1}{5}}$   \\ 
\bottomrule[1pt]
\end{tabularx}
\end{table}

To simplify notation, we denote $V_T:=V_\alpha +V_f $. 
The regret bounds of Algorithms~\ref{alg:first-order} and \ref{alg:zeroth-order} are summarized in Table~I. The first-order learning algorithm achieves a tighter regret bound than the zeroth-order counterpart for all $a>0$. 
This aligns with intuition as first-order methods leverage gradient information to guide the optimization direction and obtain more accurate updates than zeroth-order methods. However, this improvement is achieved under stronger conditions, as stated in Assumptions~\ref{assumption:J gradient bound} and \ref{assumption:probability lower bound}.

\section{Simulations}\label{sec:simulation}
This section presents a case study based on a parking lot dynamic pricing model \cite{ray2022decision}. Because parking decisions depend on both price and availability, we adjust prices dynamically in response to real-time demand. The occupancy rate $r_t \in [0,1]$ is modeled as
\begin{equation}
r_t = \xi_t + A x_t,
\end{equation}
where $x_t$ is the price, $A = -0.15$ is the price elasticity estimated from real data \cite{ray2022decision}. Let random variable $\xi_t$ follow a uniform distribution $ \mathcal{U}:=[0.9,1.1]$ to capture environmental randomness.
A time-varying target occupancy $r_t^d$ is imposed to maintain parking accessibility, which changes in response to seasonal and pedestrian traffic patterns. We select the horizon as $T =500$ and use a piecewise-constant target occupancy:
\begin{equation}\label{eq:desired occupancy step}
    r_t^d = \left\{\begin{array}{ll}
      0.65    & {\text{if}}~t \in [0,200] \\
       0.7  & {\text{if}}~ t \in (200,500].
    \end{array}\right.
\end{equation} 
The loss function is given as
\begin{equation*} 
J_t(x_t,\xi_t) = \bigl( \xi_t + A x_t - r^d_t \bigr)^2 + \frac{v}{2} x_t^2,  
\end{equation*}
with the regulation parameter $v = 0.005.$
To mitigate the risk of overcrowding, we minimize the following CVaR-based objective:
\begin{equation}\label{eq:simulation cvar}
C_t(x_t)= \mathop{\text{CVaR}_{\alpha_t}}\limits_{\xi_t \sim \mathcal{U}}[J_t(x_t,\xi_t)],
\end{equation} 
where the risk level varies consistently with the changing target occupancy rate. 
Suppose that the decision-maker becomes more conservative as the desired occupancy increases. Then, the corresponding risk level is given by 
\begin{equation}\label{eq:risk level step}
    \alpha_t = \left\{\begin{array}{ll}
      0.5    & {\text{if}}~t \in [0,200] \\
       0.8  & {\text{if}}~ t \in (200,500].
    \end{array}\right. 
\end{equation}  
The parking lot periodically adjusts its prices to monitor occupancy rates, taking into account changes in desired occupancy and risk levels. 

\begin{figure}
    \centering
\includegraphics[width=0.96\linewidth]{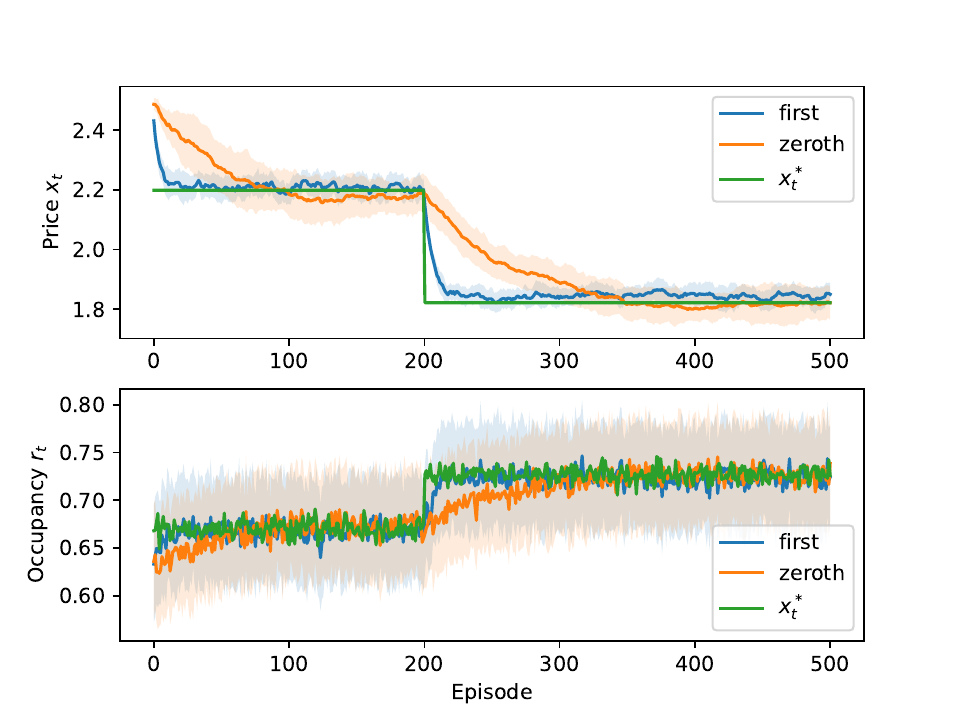}
    \caption{Top: The parking prices generated by the first-order learning algorithm (Algorithm \ref{alg:first-order}), the zeroth-order learning algorithm (Algorithm \ref{alg:zeroth-order}), and the brute-force optimal price $x_t^\ast$, under the desired occupancy \eqref{eq:desired occupancy step} and risk level \eqref{eq:risk level step}. Bottom: The resulting occupancies for each price trajectory.}
    \label{fig:price step}
\end{figure}

\begin{figure}
    \centering
    \includegraphics[width=0.96\linewidth]{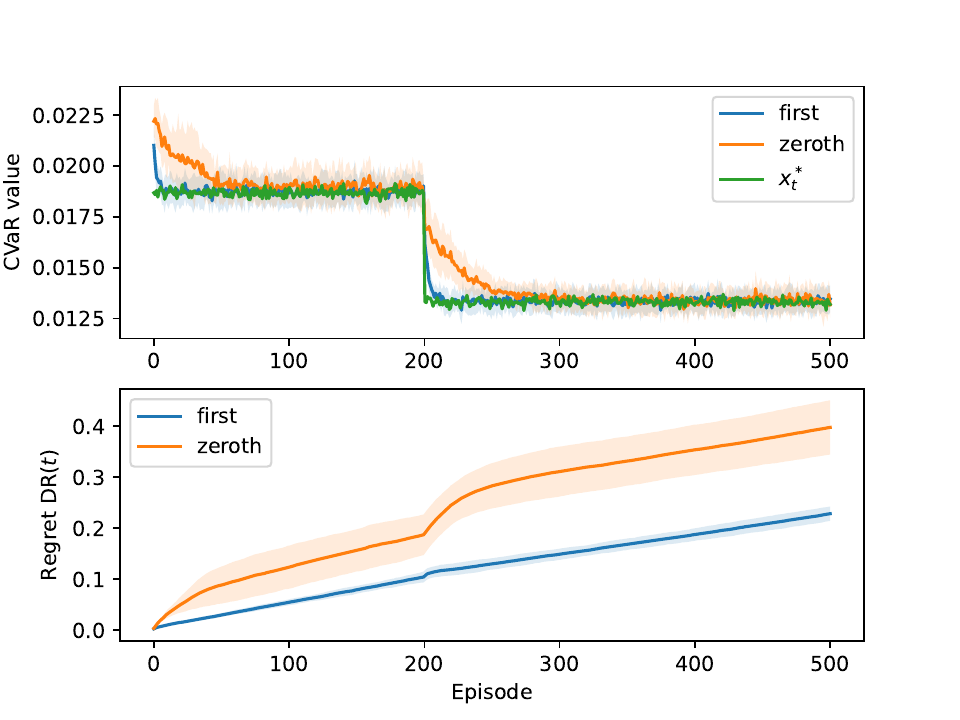}
    \caption{Top: The CVaR values under the price generated by the first-order algorithm (Algorithm~\ref{alg:first-order}),   by the zeroth-order algorithm (Algorithm~\ref{alg:zeroth-order}), and the brutal-force optimal price $x_t^\ast$, under the desired occupancy \eqref{eq:desired occupancy step} and the risk level \eqref{eq:risk level step}. Bottom: The corresponding dynamic regret.  }
    \label{fig:regret step} 
\end{figure}

Fig.~\ref{fig:price step} shows the parking prices and the corresponding occupancies generated by the first- and zeroth-order risk-averse learning algorithms (Algorithms~\ref{alg:first-order} and \ref{alg:zeroth-order}).
We set the sample number as $n_t = 8$ for all $t$. In the figures below, shaded areas represent $\pm$ one standard deviation over $20$ runs. Fig.~\ref{fig:price step} demonstrates that the prices and occupancies generated by Algorithms~\ref{alg:first-order} and \ref{alg:zeroth-order} both converge to the optimal one. Additionally,  the prices and the occupancies generated by Algorithm~\ref{alg:first-order} track desired values faster than 
Algorithm~\ref{alg:zeroth-order}.
Fig.~\ref{fig:regret step} shows the CVaR values and the dynamic regret $\mathrm{DR}(t)$ of Algorithms~\ref{alg:first-order} and \ref{alg:zeroth-order} under the changing desired occupancy \eqref{eq:desired occupancy step} and risk level \eqref{eq:risk level step}, respectively. 
The first-order risk-averse learning algorithm outperforms the zeroth-order algorithm by achieving smaller CVaR values and a lower regret.
Moreover, it can be observed from  Figs.~\ref{fig:price step}--\ref{fig:regret step} that Algorithm~\ref{alg:first-order} exhibits lower variance in its outcomes compared to Algorithm~\ref{alg:zeroth-order}, demonstrating the stabilizing effect of gradient information.

\begin{figure}
    \centering
\includegraphics[width=0.96\linewidth]{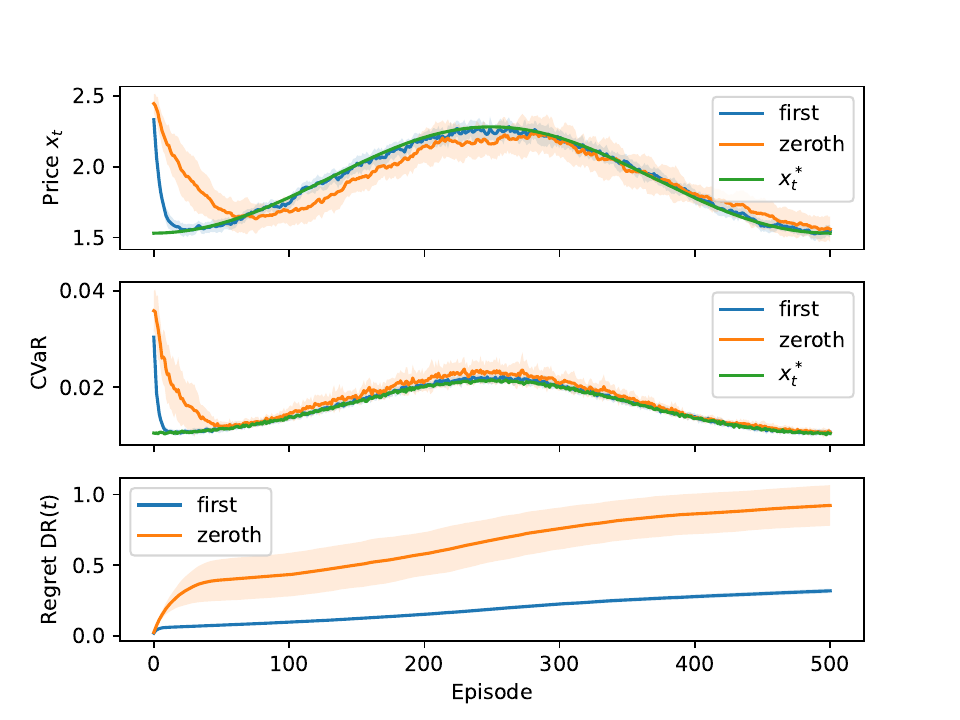}
    \caption{From top to bottom: the parking prices generated by the first-order algorithm (Algorithm~\ref{alg:first-order}), the zeroth-order algorithm (Algorithm~\ref{alg:zeroth-order}), and the brute-force optimal price $x_t^\ast$, under the target occupancy $r_t^d = 0.7 + 0.05\cos(2\pi t/500) $ and the risk level $\alpha_t = 0.5 + 0.3\cos(2\pi t/500) $; the corresponding occupancies under each price trajectory; and the resulting dynamic regret. }
    \label{fig:price sin}
\end{figure}

We next consider a case with stronger fluctuations, where the target occupancy and risk level vary sinusoidally: $r_t^d = 0.7 + 0.05\cos(2\pi t/T) $ and $\alpha_t = 0.5+0.3\cos(2\pi t/T)$.
Fig.~\ref{fig:price sin} depicts the parking prices, occupancy rates, and dynamic regret ${\rm DR}(t)$ produced by Algorithms~\ref{alg:first-order} and \ref{alg:zeroth-order}.
The results show that, under more pronounced environmental variation, the first-order algorithm still achieves faster convergence and lower regret than its zeroth-order counterpart.


\begin{figure}
    \centering
    \includegraphics[width=0.96\linewidth]{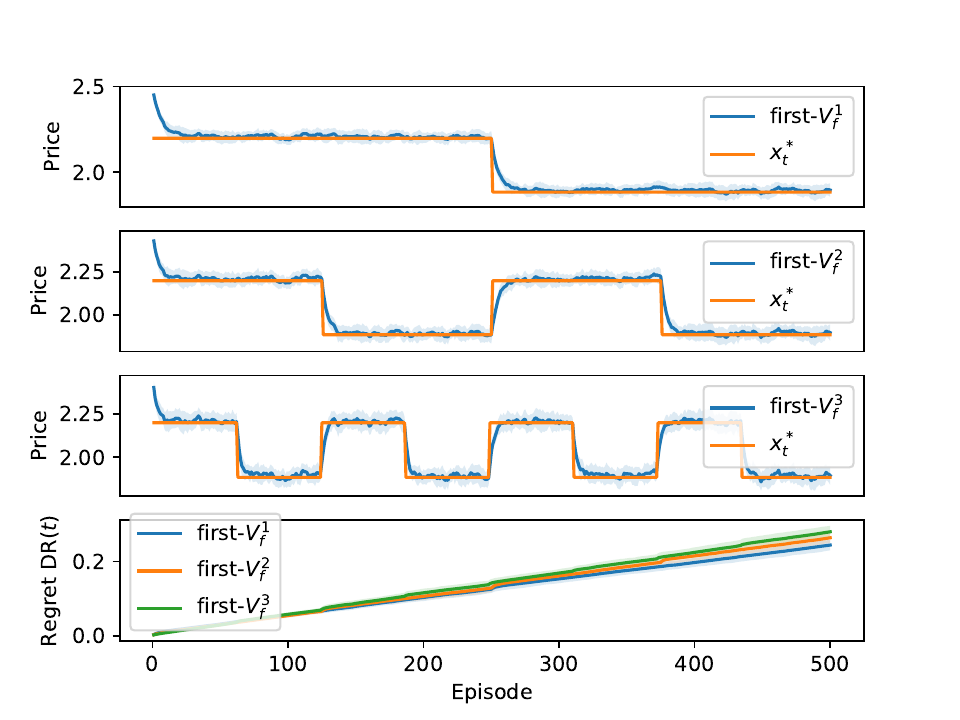}
    \caption{From top to bottom: the prices generated by Algorithm~\ref{alg:first-order} and the brutal-force optimal prices $x_t^\ast$, under the function variations $V_f^1$, $V_f^2$ and $V_f^3$, respectively; and the resulting dynamic regret. }
    \label{fig:compare Vf}
\end{figure}

To evaluate Algorithm~\ref{alg:first-order} under different function variations, we consider the desired occupancy rate as a piecewise constant function:
\begin{equation*}
r_t^d =
\begin{cases}
0.65 & \text{if } \left\lfloor \dfrac{2^m t}{T} \right\rfloor \text{ is even} \\[6pt]
0.7  & \text{otherwise},
\end{cases}
\end{equation*}
where $ m \in \{1, 2, 3\}$ controls the number of switches of $r_t^d$ within the horizon $T$, corresponding to the function variation $V_f^1$, $V_f^2$, and $V_f^3$, respectively. 
By Definition~\ref{def:function variation}, we have that $V_f^1 \le V_f^2 \le V_f^3$. 
The risk level is fixed as 
$\alpha_t=0.5$. 
For each scenario, the learning rate is tuned to achieve optimal performance, and the results are compared with the optimal price obtained via brute-force search. Fig.~\ref{fig:compare Vf} shows that Algorithm~\ref{alg:first-order} quickly tracks the optimal price under function variations $V_f^1$, $V_f^2$, and $V_f^3$.  Additionally, dynamic regret grows as the function variation increases, which aligns with Theorem~\ref{theorem:first-order}.

\begin{figure}
    \centering
\includegraphics[width=0.96\linewidth]{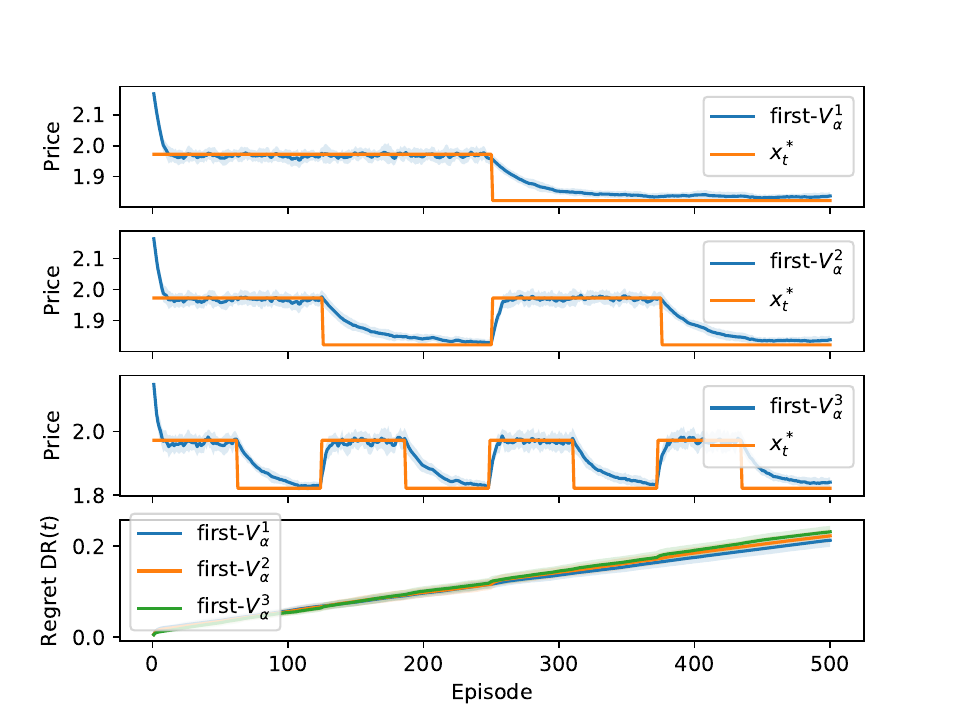}
    \caption{From top to bottom: the prices generated by Algorithm~\ref{alg:first-order} and the brutal-force optimal prices $x_t^\ast$, under the risk-level variations $V_\alpha^1$, $V_\alpha^2$ and $V_\alpha^3$,  respectively; and the resulting dynamic regret. }
    \label{fig:compare alpha}
\end{figure}

\begin{figure}
    \centering
    \includegraphics[width=0.96\linewidth]{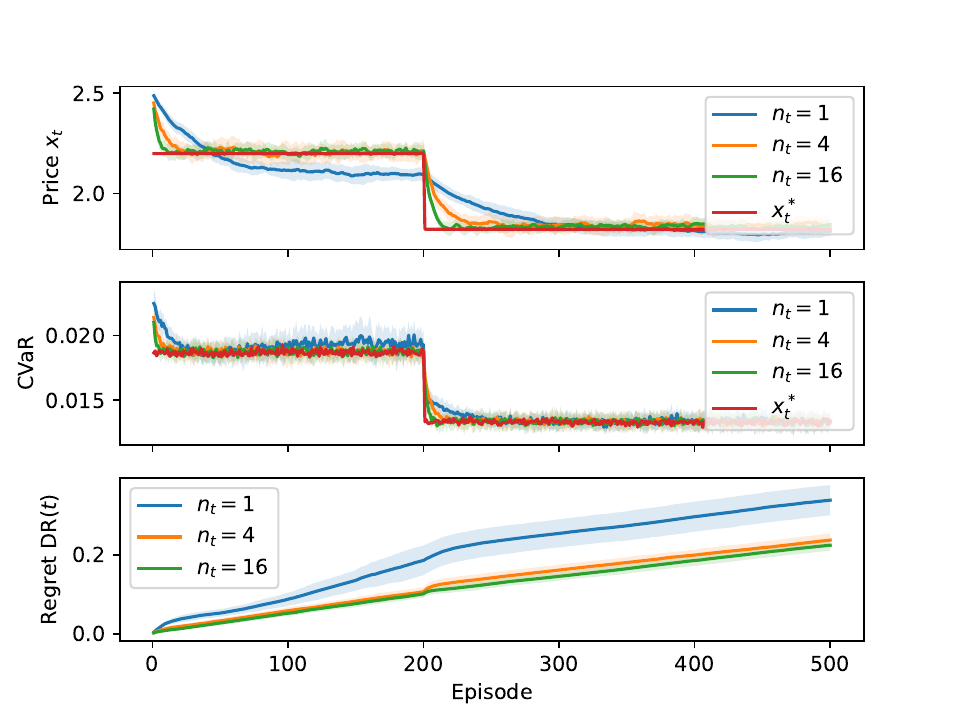}
    \caption{From top to bottom: the prices generated by Algorithm~\ref{alg:first-order} with sample strategies $n_t = \{1,4,16\}$ and the brutal-force optimal price $x_t^\ast$;  the corresponding CVaR values and the dynamic regret. }
    \label{fig:compare sample}
\end{figure}

To evaluate Algorithm~\ref{alg:first-order} under different risk-level variations, we consider the risk level as a piecewise constant function:
\begin{equation*}
\alpha_t =
\begin{cases}
0.1 & \text{if } \left\lfloor \dfrac{2^m t}{T} \right\rfloor \text{ is even} \\[6pt]
0.8  & \text{otherwise},
\end{cases}
\end{equation*}
where $ m \in \{1, 2, 3\}$ controls the number of switches of $\alpha_t$ within the horizon $T$, corresponding to $V_\alpha^1$, $V_\alpha^2$, and $V_\alpha^3$, respectively. 
By Definition~\ref{def:risk-level variation}, we have that $V_\alpha^1 \le V_\alpha^2 \le V_\alpha^3$. 
The desired occupancy is fixed as 
$r_t^d=0.7$. 
For each scenario, the learning rate is tuned to maximize performance, and the results are compared with the optimal price obtained via brute-force search. Fig.~\ref{fig:compare alpha} shows that Algorithm~\ref{alg:first-order} closely tracks the optimal price under risk-level variations $V_\alpha^1$, $V_\alpha^2$, and $V_\alpha^3$. Furthermore, dynamic regret increases with the risk-level variation, which aligns with Theorem~\ref{theorem:first-order}.
Additionally, Fig.~\ref{fig:compare sample} presents the performance of  Algorithm~\ref{alg:first-order} with sample sizes $n_t = \{1,4,16\}$ for all $t$. The results show that prices converge faster with more samples, and both CVaR values and dynamic regret decrease as the sample size increases. 

\begin{figure}
    \centering
\includegraphics[width=0.96\linewidth]{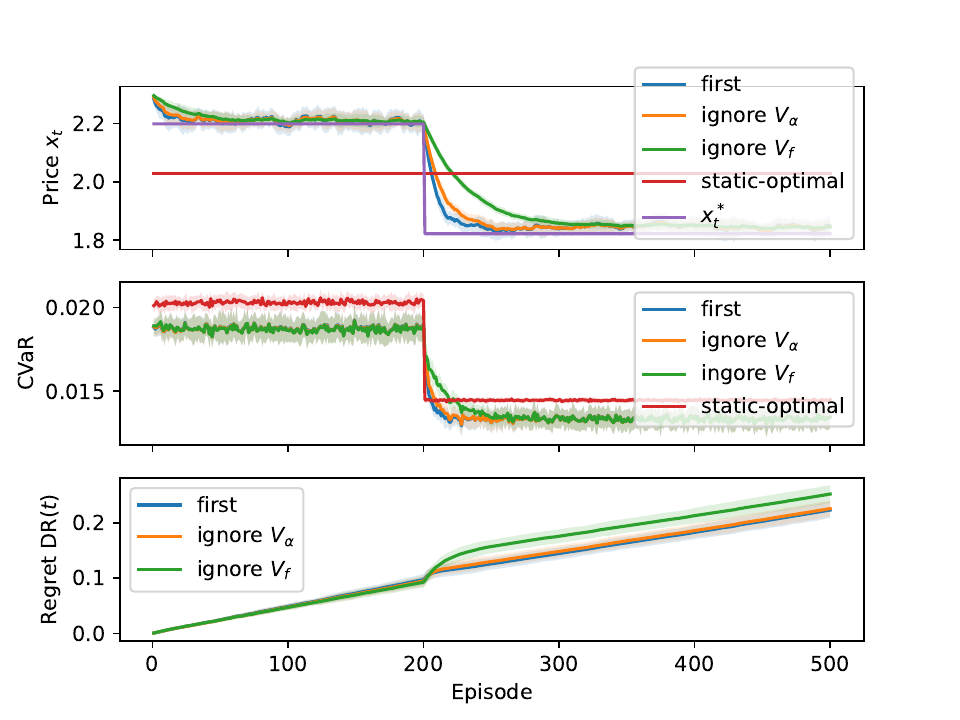}
    \caption{From top to bottom: the prices generated by Algorithm~\ref{alg:first-order}, an algorithm that ignores the risk-level variations, an algorithm that ignores the function variations, along with the optimal static price and the optimal dynamic price $x_t^\ast$ obtained via brute-force search; the corresponding CVaR values; and the resulting dynamic regret. }
    \label{fig:compare algs}
\end{figure}
Moreover, we compare Algorithm~\ref{alg:first-order} against three benchmarks: (1) a first-order learning algorithm that ignores function variation; (2) a first-order learning algorithm that ignores risk-level variation; and (3) a static scheme employing a single optimal static price.
The learning rate for each benchmark is tuned separately to achieve optimal performance. 
For reference, we include the optimal dynamic price obtained through a brute-force search. Fig.~\ref{fig:compare algs} demonstrate that Algorithm~\ref{alg:first-order} consistently outperforms all benchmarks by achieving faster convergence, lower CVaR, and reduced dynamic regret.

\section{Conclusion}\label{sec:conclusion}

This paper studied risk-averse online learning with dynamic cost functions and time-varying risk levels. 
We developed first- and zeroth-order optimization algorithms that estimate CVaR gradients from a limited number of samples. 
To capture the impact of environmental and risk dynamics, we adopted a function-variation measure and introduced a new metric for risk-level variation. 
Dynamic regret bounds for both methods are explicitly characterized in terms of function variation, risk-level variation, and the sampling budget. 
The results demonstrate the adaptability of the algorithms
in non-stationary and adaptive risk-level settings.
Moreover, the first-order algorithm achieves a tighter regret bound than the zeroth-order counterpart, confirming the benefit of gradient information.
Future work will extend this framework to risk-averse games, specifically investigating the existence of the dynamic Nash equilibria under heterogeneous and time-varying risk preferences.

\section*{Appendix}\label{sec:appendix}
\begin{lemma}
\label{lemma:function variation}\cite{besbes2015non}
Consider the function $f_t(x): \mathcal{X} \rightarrow \mathbb{R}$, define the function sequences over iteration horizon $T$ as $\{f_t(x_t)\}_{t=1}^{T}$,  we have that
\begin{align*} 
    \sum_{j=1}^s   \sum\limits_{t \in \mathcal{T}_j}  \big(f_t(\tilde{x}_j^{\ast}) -  f_t(x_t^{\ast})\big)   \leq  2  \Delta_T \sum_{t=2}^T \sup_{x\in \mathcal{X}}|f_t(x)-f_{t-1}(x)|, 
\end{align*}  
where $\tilde{x}_j^\ast = \min_{x\in\mathcal{X}}\sum_{t\in\mathcal{T}_j}f_t(x)$ and $x_t^\ast = \min_{x\in\mathcal{X}}f_t(x)$.  
\end{lemma}

\noindent\textit{Proof of Lemma~\ref{lemma:var bound}.} 
The proof of the CVaR gradient estimation error bound builds upon the approach presented in \cite{wang2024learning}. 
The concentration inequality \eqref{eq:DKW 1} follows from \cite{shao1999mathematical}. Let $\bar{\gamma} = 2e^{-2n_t\varepsilon^2\underline{p}^2}$ in \eqref{eq:DKW 1}, we have 
\begin{equation}\label{eq:var error bound}
\mathbb{P}\Bigg\{|\nu_t-\nu_t^\ast|>\frac{\sqrt{\ln{(2/\bar{\gamma})}}}{\underline{p}\sqrt{2n_t}}  \Bigg\} \le \bar{\gamma}.
\end{equation}
Denote the event in \eqref{eq:var error bound} as $A_t$, and $\mathbb{P}\{A_t\}$ the occurrence probability of $A_t$, for $t=1,\dots, T$. 

Denote $\nu_m = \min\{\nu_{t}, \nu^{\ast}_{t}\},$ and $
\nu_M = \max\{\nu_{t}, \nu^{\ast}_{t}\}$. We obtain 
\begin{align*}
e_t^1 
&= \nabla_x H_t(x_t, \nu_t) - \nabla_x C_t(x_t) \nonumber \\
&= \mathbb{E} \left[\frac{1}{\alpha_t}\mathbf{1}\{J_t(x_t, \xi) \geq \nu_t\}\nabla_x J_t(x_t, \xi)\right] \nonumber \\
&\quad - \mathbb{E}\left[\frac{1}{\alpha_t}\mathbf{1}\{J_t(x_t, \xi) \geq \nu^{\ast}_t\}\nabla_x J_t(x_t, \xi)\right]\nonumber \\
&\hspace{-1em}= \mathbb{E}\left[\frac{1}{\alpha_t}\mathrm{sgn}(\nu_t - \nu^{\ast}_t)
\mathbf{1}\{\nu_m \leq J_t(x_t, \xi) \leq \nu_M\}\nabla_x J_t(x_t, \xi)\right].
\end{align*} 
Since $\| \mathbb{E}[XY] \| \leq \mathbb{E}[\|X\|\|Y\|]$ holds for random variables $X,Y$, we have
\begin{align}\label{eq:CVaR estimate error bound}
\|e_t^1\|
&\le \frac{G}{\alpha_t} \mathbb{E} \left[\mathbf{1} \{\nu_m \le J_t(x_t, \xi) \le \nu_M\}   \right] \nonumber \\
&= \frac{G}{\alpha_t}\big(F_t(\nu_M) - F_t(\nu_m)\big)   \le \frac{G L_g}{\alpha_t} |\nu_t - \nu^{\ast}_t|  \nonumber \\
&\le \frac{GL_g\sqrt{\ln (2T/\gamma)}}{\alpha_t\underline{p}\sqrt{2n_t} }, 
\end{align}
for $t=1,\dots,T$ with probability $1-\gamma$, which establishes $    1 - \mathbb{P}\{\bigcup_{t=1}^{T} A_t \} \ge 1 - \sum_{t=1}^{T} \mathbb{P}\{A_t\} \ge 1 - T\frac{\gamma}{T} \ge 1-\gamma$. In \eqref{eq:CVaR estimate error bound}, the first inequality follows from $\|\nabla_x J_t(x_t, \xi) \|\le G$, as in Assumption~\ref{assumption:J gradient bound}. The second inequality follows from $F_t(\cdot)$ is $L_g$-Lipschitz, as in Assumption~\ref{assumption:probability lower bound}, and the last inequality follows from \eqref{eq:var error bound}. 
\hfill $\qed$

\noindent\textit{Proof of Lemma~\ref{lemma:cvar-risk variation bound}:} Define the augmented functions 
\begin{align*}
    &L_{\alpha_1}(v) = v+\frac{1}{\alpha_1} \mathbb{E}_{X\sim \mathcal{D}_x}[J(X)-v]_{+} ,\\
    &L_{\alpha_2}(v) = v+\frac{1}{\alpha_2} \mathbb{E}_{X\sim \mathcal{D}_x}[J(X)-v]_{+}.
\end{align*}
As shown in \cite{rockafellar2000optimization}, we have
\begin{align*}    \mathop{{\text{CVaR}}_{\alpha_1}}[J(X)] =  \min\limits_{v}L_{\alpha_1}(v),  \ \mathop{{\text{CVaR}}_{\alpha_2}} [J(X)] = \min\limits_{v}L_{{\alpha_2}}(v).
\end{align*}
Let $v_{\alpha_1} = \text{argmin}_{v}L_{\alpha_1}(v)$ and $v_{\alpha_2} = \text{argmin}_{v}L_{\alpha_2}(v)$.
Then, we have
\begin{equation*}
    L_{\alpha_1}(v_{\alpha_1}) = \mathop{{\text{CVaR}}_{\alpha_1}}[J(X)], \ 
     L_{\alpha_2}(v_{\alpha_2}) = \mathop{{\text{CVaR}}_{\alpha_2}}[J(X)].
\end{equation*}
Furthermore, we have
\begin{flalign*}
  & \mathop{{\text{CVaR}}_{\alpha_1}} [J(X)]  -  \mathop{{\text{CVaR}}_{\alpha_2}} [J(X)] \nonumber \\
 =  &L_{\alpha_1}(v_{\alpha_1}) -L_{\alpha_2}(v_{\alpha_2}) \nonumber \\
  \leq & L_{\alpha_1}(v_{\alpha_2}) - L_{\alpha_2}(v_{\alpha_2}) \nonumber \\ 
  = & v_{\alpha_2} + \frac{1}{{\alpha_1}}\mathbb{E}_{X\sim \mathcal{D}_x}[J(X)-v_{\alpha_2}]_{+}  \nonumber \\
  &  -v_{\alpha_2} - \frac{1}{{\alpha_2}}\mathbb{E}_{X\sim \mathcal{D}_x} [J(X)-v_{\alpha_2}]_+ \nonumber \\
 =& \Big| \frac{1}{{\alpha_1}} - \frac{1}{\alpha_2}\Big| \mathbb{E}_{X\sim \mathcal{D}_x} [J(X)-v_{\alpha_2}]_{+}
 \le  \Big| \frac{1}{{\alpha_1}} - \frac{1}{\alpha_2} \Big|  U , &&
\end{flalign*}
where the first inequality follows from $v_{\alpha_1} = \text{argmin}_{v}L_{{\alpha_1}}(v)$ and $L_{\alpha_1}(v_{\alpha_1}) \le L_{\alpha_1}(v_{\alpha_2}) $, and the second inequality from $|J(\cdot)| \le U$.
\hfill $\qed$

\noindent\textit{Proof of Lemma~\ref{lemma:cvar-function variation bound}:}
Define the augmented functions 
\begin{align*}
    &L_1(v) = v+\frac{1}{\alpha} \mathbb{E}_{X\sim \mathcal{D}_x}[J_1(X)-v]_{+} ,\\
    &L_{2}(v) = v+\frac{1}{\alpha} \mathbb{E}_{X\sim \mathcal{D}_x}[J_2(X)-v]_{+}.
\end{align*}
As shown in \cite{rockafellar2000optimization}, we have
\begin{align*}  \mathop{{\text{CVaR}}_{\alpha}}[J_1(X)] =  \min\limits_{v}L_1(v), \quad  \mathop{{\text{CVaR}}_{\alpha}} [J_2(X)] = \min\limits_{v}L_2(v).
\end{align*}
Let $v_1 = \text{argmin}_{v}L_1(v)$ and $v_2 = \text{argmin}_{v}L_2(v)$.
Then, we have
\begin{equation*}
    L_1(v_1) = \mathop{{\text{CVaR}}_\alpha}[J_1(X)], \quad  
     L_2(v_2) = \mathop{{\text{CVaR}}_\alpha}[J_2(X)].
\end{equation*}
Furthermore, we have
\begin{flalign*}
  &  \mathop{{\text{CVaR}}_{\alpha}} [J_1(X)]  -  \mathop{{\text{CVaR}}_{\alpha}} [J_2(X)] \nonumber \\
  =& L_1(v_1) -L_2(v_2)  \leq  L_1(v_2) - L_2(v_2) \nonumber \\ 
  =& v_2 + \frac{1}{{\alpha}}\mathbb{E}_{X\sim \mathcal{D}_x} [J_1(X)-v_2]_{+}  \nonumber \\ 
  & -v_2 - \frac{1}{\alpha}\mathbb{E}_{X\sim \mathcal{D}_x} [J_2(X)-v_2]_+ \nonumber \\
 = &\frac{1}{\alpha} \Big(      \mathbb{E}_{X\sim \mathcal{D}_x} [J_1(X)-v_2]_{+} -\mathbb{E}_{X\sim \mathcal{D}_x} [J_2(X)-v_2]_{+} \Big) \nonumber \\ 
 \le& \frac{1}{\alpha}        \mathbb{E}_{X\sim \mathcal{D}_x} \big[ [J_1(X)-J_2(X)]_{+}\big]   \nonumber \\ 
  \le & \frac{1}{\alpha}        \mathbb{E}_{X\sim \mathcal{D}_x} \big[ |J_1(X)-J_2(X)|\big], &&
\end{flalign*}
where the first inequality follows from $v_1 = \text{argmin}_{v}L_1(v)$ and $L_1(v_1)\le L_1(v_2)$, the second inequality  from $a_+ - b_+ \leq [a-b]_+$, and the last inequality from the 1-Lipschitz property of $[\cdot]_{+}$. A similar result appears in \cite{pichler2017quantitative}. However, we extend their analysis by relaxing the convexity assumption.
\hfill $\qed$

\bibliographystyle{ieeetr}
\bibliography{references}

\vspace{-2em}

\end{document}